\documentclass[10pt]{article}
\usepackage{verbatim}

\usepackage {   amssymb     }

\begin{document}

\title{\textbf{A new upper bound on Rubik's cube group \\ }}
\author{Silviu Radu  \\ \textbf{Lunds Institute of Technology} \\ e-mail: \texttt{silviu@bredband.net}}

\maketitle

\begin{abstract}

 In this paper a computational method is introduced for
proving that the diameter of the Rubik cube group must be less then
or equal to 40 in the quarter turn metric. Recently it has been
proven that when restricting the cube to edges the diameter is
18(see \cite{lo}). Now given an arbitrarily scrambled cube one can
solve the edges in at most 18 moves. When the edges are in place
there are 44089920 possibilities for the corners. We prove that each
of these states can be solved in at most 22 moves giving a new upper
bound of 18+22=40. The program to solve this problem is made in the
GAP language and takes about 10 hours to verify on a standard 3Ghz
PC. Later some suggestions of further improvements on the diameter
are given.
\end{abstract}
\section*{General concepts}

In this section we introduce some general notions. The  Rubik's cube
consists of 26 smaller cubies. We divide those cubies in three sets:
the edge cubies ,the corner cubies and the center cubies. The edge
cubies are those cubies that only have two colours and the corner
cubies are those that have three colours. Finally the center cubies
have only one color. So observe that an edge cubie has two faces and
a corner cubie
three faces.\\

 We will work with GAP. GAP is a computer algebra program and has a lot of special functions
 that makes handling of sets, permutations and other mathematical structures very
 easy.
 A manual with all instructions can be found at \cite{ga}. When one wants to
work with the cube group in GAP one has to specify its generators as
permutations in cycle form. One suggestion of how this can be done
is given in Martin Schoenerts example \cite{ma}. Here the same
example is used with a different
enumeration: \\
\\
$\begin{array}{cccccccccccc}
      &   &   & \mathbf{25} & \mathbf{1} & \mathbf{26} &   &   &   &   &   &   \\
      &   &   & \mathbf{8} & \mathbf{U} & \mathbf{2} &   &   &   &   &   &   \\
      &   &   & \mathbf{36} & \mathbf{4} & \mathbf{29} &   &   &   &   &   &   \\
    \mathbf{31} & \mathbf{6} & \mathbf{40} & \mathbf{27} & \mathbf{22} & \mathbf{33} & \mathbf{45} & \mathbf{17} & \mathbf{44} & \mathbf{38} & \mathbf{11} & \mathbf{46} \\
    \mathbf{3} & \mathbf{L} & \mathbf{12} & \mathbf{15} & \mathbf{F} & \mathbf{18} & \mathbf{9} & \mathbf{R} & \mathbf{13} & \mathbf{5} & \mathbf{B} & \mathbf{19} \\
    \mathbf{28} & \mathbf{20} & \mathbf{47} & \mathbf{35} & \mathbf{23} & \mathbf{43} & \mathbf{37} & \mathbf{21} & \mathbf{48} & \mathbf{30} & \mathbf{24} & \mathbf{42} \\
      &   &   & \mathbf{32} & \mathbf{16} & \mathbf{39} &   &   &   &   &   &   \\
      &   &   & \mathbf{14} & \mathbf{D} & \mathbf{10} &   &   &   &   &   &   \\
      &   &   & \mathbf{41} & \mathbf{7} & \mathbf{34} &   &   &   &   &   &
  \end{array}
$
\\

 The picture above is a picture of the cube and we have enumerated
 each corner cubie face and each edge cubie face. U stands for up, L
 for left, R for right, D for down, F for front, B for back.\\

We use the same letters for the generators as for the faces. For
example the generator U twists the face "up" clockwise. We define
the generators by typing in GAP:
\begin{verbatim}
U:=(1,2,4,8)(6,11,17,22)(25,26,29,36)(27,31,38,45)(33,40,46,44);
L:=(3,6,12,20)(8,15,14,19)(25,27,32,42)(28,31,40,47)(35,41,46,36);
F:=(4,9,16,12)(15,22,18,23)(27,33,43,35)(29,37,32,40)(36,45,39,47);
R:=(2,5,10,18)(9,17,13,21)(26,30,39,33)(29,38,34,43)(37,45,44,48);
B:=(1,3,7,13)(5,11,19,24)(25,28,34,44)(26,31,41,48)(30,38,46,42);
D:=(7,14,16,10)(20,23,21,24)(28,35,37,30)(32,39,34,41)(42,47,43,48);
\end{verbatim}

 Looking at these permutations one sees that the edges are
numerated from 1 to 24 and the corners from 25 to 48. We generate
the cube group by typing in GAP:
\begin{verbatim}
gap>cube:=Group(U,L,F,R,B,D);
 <permutation group with 6 generators>
\end{verbatim}
Now we want to let the cube act on the points 1 to 48. And we want
the orbit of this action. This is done by typing in GAP:
\begin{verbatim}
gap>Orbits(cube,[1..48],OnPoints);
[[1,2,3,4,5,6,7,8,9,10,11,12,13,14,15,16,17,18,19,20,21,22,23,24],[25,
26,27,28,29,30,31,32,33,34,35,36,37,38,39,40,41,42,43,44,45,46,47,48]]
\end{verbatim}
 The response to the last command shows that points 1 to 24(edges) are permuted among
themselves and the points 25 to 48(corners) are permuted among
themselves. This means that each group element of the cube can be
uniquely written as the product of a permutation on the set 1 to 24
and a permutation on the set 25 to 48. We will write it as
$(g_1,g_2)$ where $g_1$ is a permutation on the set 1 to 24 and
$g_2$ is a permutation on the set 25 to 48. The edge group is the
group with generators $U,L,F,B,R,D$ restricted to 1 to 24 i.e. if we
write the generators as
$(u_e,u_c)=U,(l_e,l_c)=L,(f_e,f_c)=F,(b_e,b_c)=B,(r_e,r_c)=R,(d_e,d_c)=D$
then the edge group is generated by $u_e,l_e,f_e,b_e,r_e,d_e$. We
denote it by $R_e$. The corner group is defined to be the group with
generators $u_c,l_c,f_c,b_c,r_c,d_c$. The Rubik cube group will be
denoted by $G_{cube}$. Note that neither the edge group or the
corner group is a subgroup of the cube. Actually if we for example
intersect the corner group with the cube group we get a group two
times smaller then corner group. The analogue holds for the edge
group. This can be verified in GAP. The intersection of the corner
group with the cube group we will call "the group that fixes edges"
or $E_f$ and is very important in this text.

The direct product between the corner group and the edge group gives
a group twice the size of Rubik cube group and Rubik cube group is
a subgroup of index 2 of this direct product. This can be checked with GAP.   \\

Each element in a group can be written as an word where the
generators of the group and the inverses of the generators are the
letters. The length of an element of the group is defined as the
minimal number of letters required to express it. The length of an
element $g$ is written $L(g)$ when we use notation $L(g)$ we mean
the length of an element in $G_{cube}$ with respect to generators
$\{U,L,F,B,R,D\}$. And we write $L_e(g)$ when we mean the length of
an element in the group $R_e$ with respect to generators
$\{u_e,l_e,f_e,b_e,r_e,d_e\}$.

\section*{The Symmetry Group}
Define $M$ to be the group with
generators \\
$k_1=(1,24,16,22)(2,19,10,15)(3,21,12,17)(4,11,7,23)(5,14,18,8)(6,13,20,9)
\\(25,48,32,45)(26,41,39,36)(27,44,42,37)
(28,43,40,38)(29,31,34,47)(30,35,33,46)$
\\
and \\
$k_2=(1,5,7,19)(2,10,14,8)(3,11,13,24)(4,18,16,15)(6,17,21,20)(9,23,12,22)\\
(25,44,34,28)(26,48,41,31)(27,33,43,35)(29,37,32,40)(30,42,46,38)(36,45,39,47)$.

This group acts on the set
$T=\{U,L,F,B,R,D,U^{-1},L^{-1},F^{-1},B^{-1},R^{-1},D^{-1}\}$ by
conjugation i.e. if $g\in M$ and $x\in T$ then $g^{-1}xg\in T$. The
action of
conjugation by $k_1$  can be written as a permutation with cycles:\\
$(U,B^{-1},D,F^{-1})(L,R^{-1})(F,U^{-1},B,D^{-1})(R,L^{-1})$ \\
and the action of $k_2$:\\
$(U,R,D,L)(U^{-1},R^{-1},D^{-1},L^{-1})$. \\
This action can be extended from $T$ to the cube group. It is
evident that two $M$ conjugate elements have the same length. More
details about this group can be found in \cite{ho}.

\subsection*{Equivalence relations}
There are two important equivalence relations that we will refer to often in this text:\\

\begin{enumerate}
\item Two elements $g_1,g_2$ are equivalent if there $\exists h\in M$
such that $g_1=h^{-1}g_2h$. We refer to this relation as $\sim$.

\item Two elements $g_1,g_2$ are equivalent if there $\exists h\in M$
such that $g_1=h^{-1}g_2h$
or $g_1^{-1}=h^{-1}g_2h$ . We refer to this relation as $\approx$.\\
In both cases two equivalent elements have the same length.
\end{enumerate}

\pagebreak
\section*{Calculation of the diameter of the Edge Group ($R_e$)}
We will also use that the diameter of the Cayley graph with respect
to generators $u_e,l_e,f_e,b_e,d_e,r_e$ defined above is 18. This
table is presented below and was found in \cite{lo}. The first
column gives the number of elements of lengths 0 to 18. The second
column gives the number of equivalence classes under relation $\sim$
and the second column the number of equivalence classes under
relation $\approx$.
\begin{verbatim}
 Posted to Yahoo by Tom Rokicki on  Jan 2, 2004
 Dist    Positions     Unique wrt M    Unique wrt M+inv

  0               1                1                   1
  1              12                1                   1
  2             114                5                   5
  3            1068               25                  17
  4            9819              215                 128
  5          89,392            1,886                 986
  6         807,000           16,902               8,652
  7       7,209,384          150,442              75,740
  8      63,624,107        1,326,326             665,398
  9     552,158,812       11,505,339           5,759,523
 10   4,643,963,023       96,755,918          48,408,203
 11  36,003,343,336      750,089,528         375,164,394
 12 208,075,583,101    4,334,978,635       2,167,999,621
 13 441,790,281,226    9,204,132,452       4,603,365,303
 14 277,713,627,518    5,785,844,935       2,894,003,596
 15  12,144,555,140      253,044,012         126,739,897
 16          23,716              750                 677
 17              30                3                   3
 18               1                1                   1
      ---------------    -------------    --------------
      980,995,276,800   20,437,847,376    10,222,192,146




\end{verbatim}
This calculation implies that given an arbitrary element in the cube
group we need to multiply it by at most 18 generators (in other
words to apply at most 18 moves to it) to take it in $E_f$(the group
that fixes edges).
\section*{Mapping the corner group on numbers}
It is well known that the corner group is isomorphic to
$S_8\ltimes\mathbb{Z}_3^7$. This means physically that each
configuration in the corner group can be described as a permutation
of the cubies and an orientation of the cubies. Let $\phi:R_c\mapsto
S_8\ltimes\mathbb{Z}_3^7$ be such a isomorphism. If $g\in R_c$ then
$\phi(g)$ can be uniquely written as
$\phi(g)_1\phi(g)_2],\phi(g)_1\in S_8,\phi(g)_2 \in \mathbb{Z}_3^7$.
It is obvious how to find a bijection from $Z_3^7$ to $[0,3^7-1]$.
Call this bijection $\psi_1$. There is a simple way to define a
bijective map $\psi_2:S_8\mapsto [0,8!-1]$. One can define a natural
order relation on permutations in $S_8$ which means that given two
permutations $\sigma_1,\sigma_2$ and $\sigma_1<\sigma_2$ then either
$\sigma_1(1)<\sigma_2(1)$ or if $n$ is the biggest positive integer
with the property $\sigma_1(k)=\sigma_2(k),k=1,..,n$ then
$\sigma_1(n+1)<\sigma_2(n+1)$. We simply define
$\psi_2(\sigma)=\sharp\{\tau|\tau<\sigma\}$. Now simply define a
bijection $f:R_c\mapsto [1,88179840]$ as
$f(g)=3^7\cdot\psi_2(\phi(g)_1)+\psi_1(\phi(g)_2)$.

\section*{The main objective of this paper}
In this paper we show how to construct explicit expressions in
generators of length 22 or less for each element in $E_f$("the group
that fixes edges"). This group has 44089920 elements, this can be
shown by typing in GAP: 
\begin{verbatim}
gap> cube:=Group(U,L,F,B,R,D);
<permutation group with 6 generators>
gap> corners:=Group(uc,lc,fc,bc,rc,dc);
<permutation group with 6 generators>
gap> Size(Intersection(corners,cube));
44089920
\end{verbatim}
 Since we need at most 18 moves to take an
element to $E_f$ and at most 22 moves to take an element of $E_f$ to
the identity. This means also that the diameter of the Rubik's cube group is less then or equal to 40. \\
\section*{Algorithm to express elements of $E_f$ in generators}

We need to introduce some notations. Define the following function
$p:G_{cube}\mapsto R_e$ by $p((g_1,g_2))=g_1$.We
call this map "projection on the first component".\\
 Now for
a fix $g\in R_e$ define the set $S(n,g)=\{(g,x)\in
G_{cube}|L((g,x))\leq n \}$ i.e. the set of all $\mathbf{g}$ in the
cube group that can be expressed as a product of at most $n$
generators and $p(\mathbf{g})=g$ for all $\mathbf{g}\in S(n,g)$. One
can easily see that $S(n,g)$ is empty if $L_e(g)>n$. Define
$S(k,g)\cdot S(n,h)=\{x\cdot y|x\in S(k,g),y\in S(n,h)\}$. It is
obvious that $\bigcup_{g\in R_e }S(10,g)\cdot
S(12,g^{-1})=S(22,id)$. So $\bigcup_{g\in R_e,L_e(g)\leq 10
}S(10,g)\cdot S(12,g^{-1})=S(22,id)$.\\
 Following relations are
helpful when trying to understand the algorithm we will describe:

\begin{enumerate}
\item $h^{-1}S(g,n)h=S(h^{-1}gh,n),h\in M$
\item $ h^{-1}S(k,g)\cdot S(n,g^{-1})h=S(k,h^{-1}gh)\cdot
S(n,h^{-1}g^{-1}h),h\in M$.
\item $(S(g,n))^{-1}=S(g^{-1},n)$
\end{enumerate}

In order to keep track on which elements we have found generator
expressions for, a "bit vector" is used. A "bit vector" is a vector
in which each entry can take only two values 0 or 1. GAP uses
"false" for 0 and "true" for 1. Initially all entries in this "bit
vector" are set to 1 or "true". A bijection $f:R_c\mapsto
[1,88179840]$ was defined on page 6. When generator expressions for
a given element $g\in E_f$ have been found entry number $f(g)$ is
set to 0 or "false". When we use the expression "element $g$ has
been checked of" we mean that entry $f(g)$ in the "bit vector" has
been set to "false". \\
We also need to collect information about how a given element can be
expressed in generators. In GAP we use for this a set we call
"positions". The algorithm we describe stores only one
representative from each equivalence class under relation $\approx$.
When we write "store generator expression" we mean that the
generator
expression will be saved to the set "positions". \\

Now we are ready to describe the main step of the algorithm:\\
 Define
$A$ to be the set $A=\{g\in R_e|L_e(g)\leq 10\}$. Fix a $g\in A$ and
do following steps:

\begin{enumerate}

\item Compute $S(10,g)$ and generator expressions for its elements
\item Compute $S(12,g^{-1})$ and generator expressions for its elements
\item Compute $S(10,g)\cdot S(12,g^{-1})$ and generator expressions
for its elements
\item For each element $elm$ in $S(10,g)\cdot S(12,g^{-1})$ do
following steps:
\begin{enumerate}
\item Check if $elm$ has been checked. If yes do nothing if no do
the following steps
\begin{enumerate}
\item Store generator expression of $elm$
\item Compute the orbit of $elm$ under relation $\approx$.
\item Check of every element in the orbit.
\end{enumerate}

\end{enumerate}

\item Compute $S(10,g^{-1})\cdot S(12,g)$ and generator expressions
for its elements
\item For each element $elm$ in $S(10,g^{-1})\cdot S(12,g)$ do
following steps:
\begin{enumerate}
\item Check if $elm$ has been checked. If yes do nothing if no do
the following steps
\begin{enumerate}
\item Store generator expression of $elm$.
\item Compute the orbit of $elm$ under relation $\approx$.
\item Check of every element in the orbit.
\end{enumerate}

\end{enumerate}
\end{enumerate}
Let $O(g)$ be the orbit of $g$ under relation $\approx$. Then the
above algorithm checks of every element in $\bigcup_{O(g)
}S(10,g)\cdot S(12,g^{-1})$.

The algorithm consists of repeating above procedure for enough
elements $g\in A$ which are not equivalent under relation $\approx$.
It should be noted that building $S(12,g)$ is the most time
consuming step in the algorithm. Once we have $S(12,g)$ it takes no
time to build $S(12,h^{-1}gh)$. Thats why we do the above steps for
the
whole equivalence class. \\
If we for example run the algorithm with $g$ choused as identity
than 3079007 elements are checked of.
\section*{Optimal solvers}
As mentioned above $S(22,id)$ can be written using following
formula: \[ \bigcup_{g\in R_e,L_e(g)\leq 10 }S(10,g)\cdot
S(12,g^{-1})=S(22,id)\] In the above algorithm we see how to use it
in order to compute each element in $E_f$ as a product of
generators. Obviously some elements we compute might have length
that is less then 22. To compute those we use the same algorithm
with the exception that $S(10,g)$ is replaced by $S(8,g)$. And
$S(20,id)$ can be written as: \[ \bigcup_{g\in R_e,L_e(g)\leq 8
}S(8,g)\cdot S(12,g^{-1})=S(20,id)\] Now we have all elements of
length less then or equal to 20 in $E_f$. To get only elements of
length exactly 22 we must take the difference between $S(22,id)$ and
$S(20,id)$. We continue computing the set $S(18,id)$ using the
formula:
\[ \bigcup_{g\in R_e,L_e(g)\leq 8 }S(8,g)\cdot
S(10,g^{-1})=S(18,id)\]\\
And finally we compute $S(16,id)$, $S(14,id)$, $S(12,id)$, $
S(10,id)$ and $S(8,id)$. To compute elements of length $n$ we simply
take the difference $S(n,g)-S(n-1,g)$. The program for this
calculation is not given here but we give the result on next page.
Observe that 22q means length 22.
\begin{verbatim}

Distance        Nr of pos   Unique wrt M    Unique wrt M + inv

 0q             1           1               1
 2q             0           0               0
 4q             0           0               0
 6q             0           0               0
 8q             240         5               3
 10q            288         6               3
 12q            8764        197             113
 14q            116608      2475            1303
 16q            840513      17620           9091
 18q            9884342     206424          104523
 20q            30623418    639362          323859
 22q            2615746     54895           28528
                --------    -----           ------
                44089920    920980          467424
\end{verbatim}

We have also computed the length of all elements in the "Group that
fixes cubies". This is the set of all elements in the Rubik cube
group such that the cubies are in their correct spot but their
orientation is arbitrary. This is an abelian normal subgroup of the
cube. This can be checked with GAP. The order of this group is
4478976. Here is the table with lengths of each element:
\begin{verbatim}
Distance        Nr of pos  Unique wrt M  Unique wrt M + inv

 0q             1          1             1
 2q             0          0             0
 4q             0          0             0
 6q             0          0             0
 8q             0          0             0
 10q            0          0             0
 12q            441        11            8
 14q            3944       87            52
 16q            32110      708           396
 18q            456025     9656          5009
 20q            2873194    60399         30978
 22q            1113236    23652         12424
 24q            25         4             4
                --------   -----         ------
                4478976    94518         48872
\end{verbatim}
One of the positions of length 24 is the well known "superflip".
This position alone generates the center of the Rubik cube group.
The other positions up to conjugation with M are:
\begin{verbatim}
 F  R' U  B2 U  R  D  F  U' B  R  B2 D  B  R' U' F' B  R  U2 D' (24q*,21f)
 F  R' F2 R' U  F  U  L  B2 R  B  D2 F  R' B' D  B  L' D' R' U' (24q*, 21f)
 F  R  D2 B' R' F  R' D2 B' D' L  F  L2 U  B' U  F L2 U' B'     (24q*, 20f)
\end{verbatim}
This positions were the only ones not computable in 22q and the
expressions in generators were found using Michael Reids solver see
\cite {re}.
\section*{Upperbounds on different subsets}
Take an element $g\in R_e$ and $L(g)=n$  ($L(g)$ is the length of
$g$) in $R_e$. Then for any element $(g,y)\in G_{cube}$ we need at
most $n$ moves to make it an element of the form $(id,x)\in
G_{cube}$.And an element $(id,x)\in G_{cube}$ needs at most 22 moves
applied on it to make it the identity. So $L((g,y))\leq 22+n$ in
$G_{cube}$.There is only one element $h$ in $R_e$ with $L(h)=18$ in
$R_e$(according to the above diagram). So the only elements in
$G_{cube}$ that may need 40 moves are those in the set $S=\{(h,x)\in
G_{cube}\}$.If a solver solves each position in $S$ in less then 40
moves then one gets a new upper bound 39. Actually one only needs to
solve a position from each equivalence class under relation
$\approx$ which gives about 500000 positions needed to be solved.
This computation has not been done here. But soon after the author
posted generator expressions for elements in $E_f$ this has been
done by another author. See \cite {fo} for details.
\section*{The program}
To begin with we are going to create some sets we are going to work
with. The sets will be presented in GAP as lists. The first one is $\mathbf{s}$. \\
 This set is made by typing in gap:
\begin{verbatim}
 gap>s:=[U,L,F,B,R,D,U^-1,L^-1,F^-1,B^-1,R^-1,D^-1];
\end{verbatim}
The other sets are: \\
\subsection*{Some sets}
$\mathbf{veven}$ the set of all elements in the edge group of length
0, 2, 4, 6.\\
$\mathbf{vodd}$ the set of all elements in the edge group of
length 1, 3, 5. \\
$\mathbf{vevenr}$ the set of all elements in the Rubik cube group
of length 0, 2, 4, 6. \\
$\mathbf{voddr}$ the set of all elements in the Rubik cube group
of length 1, 3, 5. \\
The elements of this sets are stored as permutations with disjoint
cycles.

$\mathbf{vevenc}$ Contains at entry $k$ the same element that
$\mathbf{vevenr}$ contains at entry $k$ but stored as products of
generators. Throughout this program we store a product of generators
as a vector whose entries are integers between 1 and 12. Integer $n$
means generator $s[n]$ and $\mathbf{s}$ was defined above.\\
$\mathbf{voddc}$ is analogue to $\mathbf{vevenc}$ but the set
$\mathbf{vevenr}$ is replaced with $\mathbf{voddr}$.

 All the sets containing permutations with cycles are sorted with respect to the
natural order relation on permutations which means that given two
permutations $\sigma_1,\sigma_2$ and $\sigma_1<\sigma_2$ then either
$\sigma_1(1)<\sigma_2(1)$ or if $n$ is the biggest positive integer
with the property $\sigma_1(k)=\sigma_2(k),k=1,..,n$ then
$\sigma_1(n+1)<\sigma_2(n+1)$. \\
The elements in the set $\mathbf{vevenr}$  are ordered in following
way:\\
$(g_1,h_{k_1}),(g_1,h_{k_2}),..,(g_1,h_{k_m}),(g_2,h_{l_1}),(g_2,h_{l_2}),...$.
This is because the first component is a permutation on the set 1 to
24 while the second component a permutation on the set 25 to 48. Then the set $\mathbf{veven}$ is ordered in following way: \\
$g_1,g_2,..$. Note that we refer to an element of $\mathbf{vevenr}$ as an element of the form $(g,h)$ although it is stored as
a permutation with disjoint cycles.  \\
We use two more vectors $\mathbf{y1}$ and $\mathbf{y2}$ where the
first has the same size as $\mathbf{vodd}$ and the second the same
size as $\mathbf{veven}$.The vector $\mathbf{y1}$ contains at
position $k$ a number that gives the position in $\mathbf{voddr}$
where we have the first occurrence of an element of the form
$(\mathbf{vodd}(k),h)$.\\
Finally we call the "bit vector" for $\mathbf{rcv}$ which has the
same length as $R_c$ (88179840 elements). When the program starts
all its entries are set to true.
\subsection*{Some functions}
We will predefine some functions in order to make the program code
easier:\\

$\mathbf{MakeSetv(n)}$ if n=1 then it returns set
$\mathbf{[vevenr,voddr]}$ if n=2
returns set $\mathbf{[veven,vodd]}$. \\
$\mathbf{MakeSetvc()}$ returns $\mathbf{[vevenc,voddc]}$.\\
$\mathbf{MakeSety(n)}$ if n=1 returns set $\mathbf{y1}$ if n=2
returns set $\mathbf{y2}$.\\
 $\mathbf{MakeSet(n,g)}$ This function will give a set of
the form $(g,h), h\in R_c$ and all its elements of length $n$ or
less. The returned set is of the form $\mathbf{[set1,set2]}$ where
the elements of $set1$ are expressed as the elements of
$\mathbf{vevenr}$ and the elements of
$set2$ are expressed as products of generators.\\
$\mathbf{PermToNr(g)}$ This is a bijective function
$PermToNr:R_c\mapsto [1,88179840]$.\\
$\mathbf{EqClassRel2(el)}$ Returns the set $\{PermToNr(g)|g\approx el\}$. \\
$\mathbf{WriteRcv(set1)}$ For each $x\in set1$ does
$rcv(x)=false$.\\
 $\mathbf{Redf(set1)}$.Partitions set
$set1$ into equivalence classes under relation $\approx$ defined
above and
returns a set containing a representative from each class. \\
$\mathbf{CharToPerm(x)}$ is a function that converts a text string
to a permutation. This was done so that the attached data don't take
to much space in this paper. \\
$\mathbf{MultEl(x)}$ if for example $x:=[1,2,1,4,11,2]$ then this
function returns $s[1]*s[2]*s[1]*s[4]*s[11]*s[2]$.\\
$\mathbf{Invg(x)}$ if for example $x:=[1,2,1,4,11,2]$ then it
returns $[8,5,10,7,8,7]$ i.e.
$\mathbf{MultEl(x)^{-1}=MultEl(Invg(x))}$.
 \pagebreak
\section*{The Code}
Before we describe the code some GAP commandos must be explained.
Those can also be found in the manual \cite{ga} but we describe them
here for completeness. To build a set in GAP one types:
\begin{verbatim}
gap>set:=[];
\end{verbatim}
To build a set of sets type:
\begin{verbatim}
gap>set:=[[],[],[],[]];
\end{verbatim}
The above set contains 4 empty subsets. To store the identity in the
first subset type in GAP:
\begin{verbatim}
gap>set[1][1]:=();
\end{verbatim}
Now the first set in the set "set" contains the identity.\\
 To build
a "bit vector" of size 1000 with all entries set to "true" or 1 just
type:
\begin{verbatim}
gap>vec:=BlistList([1..1000],[1..1000]);
\end{verbatim}
To build a "bit vector" with all entries false just type:
\begin{verbatim}
gap>vec:=BlistList([1..1000],[]);
\end{verbatim}
If we type in GAP:
\begin{verbatim}
 gap>set:=[1,1];
\end{verbatim}
then we have build a set that contains the same element on two
different entries, maybe it is better to call it a list. But if we
type:
\begin{verbatim}
 gap>set2:=Union(set,set);
\end{verbatim}
 then the set "set2" will contain only element 1. If we instead of
union would have taken intersection the same thing should have
happened. A union, intersection or difference of two sets will never
return a set with doubles. Furthermore after performing those
operations the set will even be sorted.

\subsection*{Generators and the Symmetry Group}
We define all generators we are going to use and the symmetry group
in following file:
\begin{verbatim}
#Generators for the cube group
U:=(1,2,4,8)(6,11,17,22)(25,26,29,36)(27,31,38,45)(33,40,46,44);
L:=(3,6,12,20)(8,15,14,19)(25,27,32,42)(28,31,40,47)(35,41,46,36);
F:=(4,9,16,12)(15,22,18,23)(27,33,43,35)(29,37,32,40)(36,45,39,47);
R:=(2,5,10,18)(9,17,13,21)(26,30,39,33)(29,38,34,43)(37,45,44,48);
B:=(1,3,7,13)(5,11,19,24)(25,28,34,44)(26,31,41,48)(30,38,46,42);
D:=(7,14,16,10)(20,23,21,24)(28,35,37,30)(32,39,34,41)(42,47,43,48);
#Generators for edge group an corner group
ue:=RestrictedPerm(U,[1..24]);uc:=RestrictedPerm(U,[25..48]);
le:=RestrictedPerm(L,[1..24]);lc:=RestrictedPerm(L,[25..48]);
fe:=RestrictedPerm(F,[1..24]);fc:=RestrictedPerm(F,[25..48]);
re:=RestrictedPerm(R,[1..24]);rc:=RestrictedPerm(R,[25..48]);
be:=RestrictedPerm(B,[1..24]);bc:=RestrictedPerm(B,[25..48]);
de:=RestrictedPerm(D,[1..24]);dc:=RestrictedPerm(D,[25..48]);
#Generators for the group of symmetries
k1:=(1,24,16,22)(2,19,10,15)(3,21,12,17)(4,11,7,23)(5,14,18,8)(6,13,20,9)
(25,48,32,45)(26,41,39,36)(27,44,42,37)(28,43,40,38)(29,31,34,47)(30,35,33,46);
k2:=(1,5,7,19)(2,10,14,8)(3,11,13,24)(4,18,16,15)(6,17,21,20)(9,23,12,22)
(25,44,34,28)(26,48,41,31)(27,33,43,35)(29,37,32,40)(30,42,46,38)(36,45,39,47);
#Generators for the group of symmetries restricted to edge resp corners.
k1e:=RestrictedPerm(k1,[1..24]);k2e:=RestrictedPerm(k2,[1..24]);
k1c:=RestrictedPerm(k1,[25..48]);k2c:=RestrictedPerm(k2,[25..48]);
#symmetry group 
M:=Group(k1,k2);
#symmetry group restricted to corners
Mc:=Group(k1c,k2c);
#symmetry group restricted to edges
Me:=Group(k1e,k2e);
\end{verbatim} This code will be saved to file gens.txt.
 \pagebreak
\subsection*{Function MakeSetv}
This function builds the sets veven, vodd, vevenr, voddr. The
function first selects generators depending on if we want to build
veven, vodd or vevenr, voddr. If we want to build vevenr, voddr we
type:
\begin{verbatim}
gap>a:=MakeSetv(1);
\end{verbatim}
The function will return a set of sets containing the sets vevenr
and voddr. In the above case a[1] is vevenr and a[2] is voddr. Lets
describe the algorithm that this program uses:
\begin{enumerate}
\item First a set of sets is build containing 6 empty sets and a set containing the identity:
\begin{verbatim}
v:=[[()],[],[],[],[],[],[],[]]
\end{verbatim}

\item Now take each element v[1] and multiply it be each
generator and store the result in the next set i.e. v[2].
\item Now take v[2] and multiply each element of it with the
generators and store the result in v[3]. Remove doubles from v[3]
and elements contained in v[1].
 \item Now take the next set in v i.e. v[3] multiply each element
 by generators, store the result in v[4] remove doubles from v[4] and elements
 contained in v[2].
 \item Continue doing this procedure until v[6] i.e. v[6] is the
 last set of which elements will me multiplied by generators and the
 result stored in v[7]
 \item Now return the set
 \{Union(v[1],v[3],v[5],v[7]),Union(v[2],v[4],v[6])\} i.e.
 \{veven,vodd\}.

\end{enumerate}
See next page for the code of above described function:

\pagebreak

\begin{verbatim}
MakeSetv:=function(n)
local veven,vodd,vevenr,voddr,gens,selct,v,k,x,z,y,cnt;
gens:=[[()],[()]];selct:=n;
gens[1]:=[U,L,F,B,R,D,U^-1,L^-1,F^-1,B^-1,R^-1,D^-1];
gens[2]:=[ue,le,fe,be,re,de,ue^-1,le^-1,fe^-1,be^-1,re^-1,de^-1];
v:=[[()],[()],[()],[()],[()],[()],[()]];
for y in [1..6] do
    k:=[[()],[()],[()],[()],[()],[()],[()],[()],[()],[()],[()],[()]];
    cnt:=1;
    for x in v[y] do
         for z in [1..12] do
             k[z][cnt]:=x*gens[selct][z];
         od;
         cnt:=cnt+1;
    od;
    v[y+1]:=Union(k[1],k[2],k[3],k[4],k[5],k[6],k[7],k[8],k[9],k[10],k[11],k[12]);
    if y>1 then v[y+1]:=Difference(v[y+1],v[y-1]);fi;
od;
if selct=1 then 
    vevenr:=Union(v[1],v[3],v[5],v[7]);voddr:=Union(v[2],v[4],v[6]);
    return [vevenr,voddr];
elif selct=2 then 
    veven:=Union(v[1],v[3],v[5],v[7]); vodd:=Union(v[2],v[4],v[6]);
    return [veven,vodd];
fi;
end;;
\end{verbatim} This code will be saved to file
makesetv.txt.
 \pagebreak
\subsection*{Function MakeSetvc}
This function will build the sets vevenc and voddc. We remind that
this sets contain the same number of elements as vevenr and voddr.
We have defined a list "s" above as
$s=[U,L,F,B,R,D,U^{-1},L^{-1},F^{-1},B^{-1},R^{-1},D^{-1}]$. To
define it in GAP just type:
\begin{verbatim}
gap>s:=[U,L,F,B,R,D,U^-1,L^-1,F^-1,B^-1,R^-1,D^-1];
\end{verbatim}
We now want to store generator expressions for elements in vevenr
and voddr in sets vevenc and voddc if for instance vevenc contains
at position 10 generator expression [1,2,7] then vevenr should
contain at position 10 element s[1]*s[2]*s[7]. We will describe an
algorithm how to achieve this. But first we need a few preparations.
 GAP has a function called "Tuples" that we will use. We explain it
 by giving an example:
 \begin{verbatim}
gap>Tuples([1..5],2);
 \end{verbatim}
After typing this GAP will return a set containing all possible
elements of form $(x,y)$ with $x\in\{1,2,3,4,5\},y\in\{1,2,3,4,5\}$.
And we also use a function called PositionSet. We explain it using
an example:
\begin{verbatim}
gap>PositionSet(veven,(1,2,3));
\end{verbatim}
If (1,2,3) is contained in veven then this function will return a
number $n$ such that veven[n]=(1,2,3). If it is not contained it
will just return fail.  Now we are ready to describe the algorithm.
It is described only for vevenc since voddc is analogue:
\begin{enumerate}
\item Build a "bit vector" with same size a vevenr with all elements
set to "false".
\item Build the set of all tuples on elements
\{1,2,3,4,5,6,7,8,9,10,11,12\} and length 2. Call this set Tu2.
\item Take each element $(x,y)\in Tu2$ and find the position in vevenr of the element
$s[x]*s[y]$. Call this position n. If entry number n is "false" in
the "bit vector", set it to "true" and put element $(x,y)$ in vevenc
at position n. If entry number n in the "bit vector" is "true" then
do nothing.
\item Repeat step 2 and 3 for tuples of length 4 and 6 on the set
\{1,2,3,4,5,6,7,8,9,10,11,12\}.
\end{enumerate}
The code for this function is given on next page.
\begin{verbatim}
s:=[U,L,F,B,R,D,U^-1,L^-1,F^-1,B^-1,R^-1,D^-1];
MakeSetvc:=function()
local binv,binv2,vevenc,voddc,x,set,z,elm,g,pos;
binv:=BlistList([1..Size(vevenr)],[]);
binv2:=BlistList([1..Size(voddr)],[]);
vevenc:=[];
voddc:=[];
for x in [1..6] do
if x=1 or x=3 or x=5 then
	set:=Tuples([1..12],x);
	for z in set do
		elm:=(); 
		for g in [1..Size(z)] do
           	elm:=elm*s[z[g]];
		od;
	pos:=PositionSet(voddr,elm);
	if binv2[pos]=false then voddc[pos]:=z;binv2[pos]:=true;fi;
	od;
fi;


if x=2 or x=4 or x=6 then
	set:=Tuples([1..12],x);
	for z in set do
	elm:=();
	for g in [1..Size(z)] do
	elm:=elm*s[z[g]];
	od;
	pos:=PositionSet(vevenr,elm);
	if binv[pos]=false then vevenc[pos]:=z;binv[pos]:=true;fi;
	od;
fi;
od;
vevenc[1]:=[];
return[vevenc,voddc];
end;;
\end{verbatim}
 \pagebreak
\subsection*{Function MakeSety}

As described before vector y1 has the same size as vodd and vector
y2 the same size as veven. We only describe how to build y1 since y2
is analogue. Remember that if an element $g\in R_e$ has position $n$
in veven then y[n] contains a number $m$ that gives the position in
vevenr of the first element who's projection on $R_e$ is $g$. By
projection we mean the map $(g,x)\mapsto g$. We use following
algorithm to build y2:
\begin{enumerate}
\item Build a vector "vtemp" the same size as vevenr that contains at
position $n$ the projection of vevenr[n] on $R_e$.
\item Do following for each $g\in veven$.
\begin{enumerate}
\item Find position of $g$ in veven call it "pos"
\item Find the position of the first occurrence of $g$ in vtemp call
it "pos2".
\item Set y2[pos]=pos2
\end{enumerate}
\end{enumerate}
The code is described below:

\begin{verbatim}
MakeSety:=function(n) local vtemp,v,vr,x,y,cnt,g;
vtemp:=[()];cnt:=1;
if n=1 then v:=vodd;vr:=voddr;
elif n=2 then v:=veven;vr:=vevenr;
fi;
x:=[1..Size(vr)];y:=[1..Size(v)];
for g in x do
   vtemp[g]:=RestrictedPerm(vr[g],[1..24]);
od;
for g in v do
   y[PositionSorted(v,g)]:=PositionSorted(vtemp,g);
od;
return y; end;;
\end{verbatim} This code will be saved to file
makesety.txt
\subsection*{Function PermToNr}
Here is the program described in the section "How to map a
permutation on a number". It takes a $g\in R_c$ and maps it on
$S_8\ltimes \mathbb{Z}_3^7$. Lets describe how to map a permutation
of $S_8$ on numbers [1,8!] in more detail. Given a $\sigma\in S_8$
we want to determine all elements that are less then $\sigma$ with
respect to the order relation that was defined on page 11. This is a
simple combinatorial problem and we give an example of how to
determine the number of elements less then a given permutation and
then generalize this to an arbitrary one. Take for example $\sigma=
\left(
                                                          \begin{array}{cccccccc}
                                                            1 & 2 & 3 & 4 & 5 & 6 & 7 & 8 \\
                                                            5 & 6 & 3 & 4 & 2 & 1 & 8 & 7 \\
                                                          \end{array}
                                                        \right)
$. We divide the set $P$ of permutations less than $\sigma$ in
following disjoint subsets:
\begin{enumerate}
\item The subset with
$\{{\tau}(1)<5,\tau\in P\}$. Then ${\tau}(1)$ can be any element in
the set $\{1,2,3,4\}$ so this set has size $7!\cdot 4$
\item The subset
with $\{{\tau}(1)=5,\tau \in P\}$. Then ${\tau}(2)$ can be any
element in the set $\{1,2,3,4\}$ so the size of this set is $6!\cdot
4$.
\item The subset with $\{{\tau}(1)=5,{\tau}(2)=6,\tau\in P\}$
then ${\tau}(3)$ can be any element in the set $\{1,2\}$ and the
size of this set is $5!\cdot 2$.
\item The subset with $\{{\tau}(1)=5,{\tau}(2)=6,{\tau}(3)=3,\tau\in
P\}$ then ${\tau}(4)$ can be any element in the set $\{2,1\}$. The
size of this subset is $4!\cdot 2$.
\item The subset with
$\{{\tau}(1)=5,{\tau}(2)=6,{\tau}(3)=3,{\tau}(4)=4\}$ then
${\tau}(5), \tau \in P$ is in the set $\{1\}$. The size of this set
is $3!$.
\item The subset with
$\{{\tau}(1)=5,{\tau}(2)=6,{\tau}(3)=3,{\tau}(4)=4,{\tau}(5)=2,\tau
\in P\}$ and ${\tau}(6)$ can be any element in the set $\{\}$. The
size of this subset is 0.
\item The subset with
$\{{\tau}(1)=5,{\tau}(2)=6,{\tau}(3)=3,{\tau}(4)=4,{\tau}(5)=2,{\tau}(6)=1,\tau\in
P\}$ then ${\tau}(7)$ can be any element in the set $\{7\}$. The
size of this subset is $1!$.
\end{enumerate}
Since all of this subsets are disjoint and their union is $P$ it
follows that the size of $P$ is $4\cdot 7!+4\cdot 6!+2\cdot
5!+2\cdot 4!+ 3!+0 + 1!$. This can be easily generalized for an
arbitrary $\sigma$. Let $\sigma=\left(
                                  \begin{array}{cccccccc}
                                    1 & 2 & 3 & 4 & 5 & 6 & 7 & 8 \\
                                    \sigma(1) & {\sigma}(2) & {\sigma}(3) & {\sigma}(4) & {\sigma}(5) & {\sigma}(6) & {\sigma}(7) & {\sigma}(8) \\
                                  \end{array}
                                \right)$ then the number of
                                permutations less than $\sigma$ is

                                \[\sum_{k=1}^7 (7-k+1)!\cdot ({\sigma}(k)-1)\cdot \sharp\{ {\sigma}(1)<{\sigma}(k)\mid l\in\{1,..,k-1\}
                                \}\]
The code for bijection $f:R_c\mapsto [1,88179840]$ can be seen on
next page.
 \pagebreak
\begin{verbatim}
PermToNr:=function(g) local m1,m2,m3,m4,m5,m6,m7,m8,y,w,nr,k1,k2;
y:=[1..48];
y[25]:=1;y[31]:=2;y[46]:=3;y[26]:=1;y[38]:=2;y[44]:=3;y[27]:=1;y[36]:=2;y[40]:=3;
y[28]:=1;y[41]:=2;y[42]:=3;y[29]:=1;y[33]:=2;y[45]:=3;y[30]:=1;y[34]:=2;y[48]:=3;
y[32]:=1;y[35]:=2;y[47]:=3;y[37]:=1;y[39]:=2;y[43]:=3; w:=[1..48];
w[25]:=1;w[31]:=1;w[46]:=1;w[26]:=2;w[38]:=2;w[44]:=2;w[27]:=3;w[36]:=3;w[40]:=3;
w[28]:=4;w[41]:=4;w[42]:=4;w[29]:=5;w[33]:=5;w[45]:=5;w[30]:=6;w[34]:=6;w[48]:=6;
w[32]:=7;w[35]:=7;w[47]:=7;w[37]:=8;w[39]:=8;w[43]:=8;
m1:=w[25^g];m2:=w[26^g];m3:=w[27^g];m4:=w[28^g];m5:=w[29^g];m6:=w[30^g];m7:=w[32^g];
m8:=w[37^g]; if m1<m2 then m2:=m2-1; fi; if m1<m3 then m3:=m3-1; fi;
if m2<m3 then m3:=m3-1; fi; if m1<m4 then m4:=m4-1; fi; if m2<m4
then m4:=m4-1; fi; if m3<m4 then m4:=m4-1; fi; if m1<m5 then
m5:=m5-1; fi; if m2<m5 then m5:=m5-1; fi; if m3<m5 then m5:=m5-1;
fi; if m4<m5 then m5:=m5-1; fi; if m1<m6 then m6:=m6-1; fi; if m2<m6
then m6:=m6-1; fi; if m3<m6 then m6:=m6-1; fi; if m4<m6 then
m6:=m6-1; fi; if m5<m6 then m6:=m6-1; fi; if m1<m7 then m7:=m7-1;
fi; if m2<m7 then m7:=m7-1; fi; if m3<m7 then m7:=m7-1; fi; if m4<m7
then m7:=m7-1; fi; if m5<m7 then m7:=m7-1; fi; if m6<m7 then
m7:=m7-1; fi;
k1:=(m1-1)*7*6*5*4*3*2+(m2-1)*6*5*4*3*2+(m3-1)*5*4*3*2+(m4-1)*4*3*2+
(m5-1)*3*2+(m6-1)*2+m7-1;
k2:=(y[25^g]-1)*3^0+(y[26^g]-1)*3^1+(y[27^g]-1)*3^2+(y[28^g]-1)*3^3+
(y[29^g]-1)*3^4+(y[30 ^g]-1)*3^5+(y[32^g]-1)*3^6; nr:=2187*k1+k2+1;
return nr; end;;
\end{verbatim} This code will be saved to file
permtonr.txt.\\
 Nobody is expected to follow this code. That is why a program is given to verify that this function is a bijection
and it is all that is needed for the proof:

\begin{verbatim}
corners:=Group(uc,lc,fc,bc,rc,dc); rcv:=BlistList([1..88179840],[]);
for x in corners do
 rcv[PermToNr(x)]:=true;
od;
Print(rcv=BlistList([1..88179840],[1..88179840]),"\n");

\end{verbatim}

If the program returns true then it is a surjection and consequently
a bijection.
 \pagebreak
\subsection*{The function MakeSet(n,g)}
This function builds $S(10,g)$ and $S(12,g)$. To build $S(10,g)$
just type:
\begin{verbatim}
gap>MakeSet(10,g);
\end{verbatim}
We explain the algorithm for $S(12,g)$
 \begin{enumerate}
\item Build an empty set $S$
\item For each element $h$ in veven do
\begin{enumerate}
\item Build the set $S_h=\{(h,x)\in vevenr\}$ and the set
$S_{h^{-1}g}=\{(h^{-1}g,x)\in vevenr\}$.
\item Build the set $\{g_1\cdot g_2\mid g_1\in S_h,g_2\in
S_{h^{-1}g}\}$ and put each element from this set in $S$.
\end{enumerate}
\item Return the set $S$ (or the set $S=S(12,g)$
 \end{enumerate}
The implementation of this function as code can be seen on next
page. \pagebreak
\begin{verbatim}
MakeSet:=function(n,g) local set,y,veo,veor,a,b,z,x,kl,j,k,l,m,t,set2,set3,nre,veoc;
if n=10 then 
 y:=y1;veo:=vodd;veor:=voddr;veoc:=voddc;
elif n=12 then
 y:=y2;veo:=veven;veor:=vevenr;veoc:=vevenc;
fi;
set:=[];t:=[1..Size(veo)];kl:=Size(veo);set2:=[];set3:=[];nre:=0;
for x in t do
  if x=kl then
    a:=Size(veor)+1-y[x];
  else
    a:=y[x+1]-y[x];
  fi;
  z:=PositionSet(veo,veo[x]^-1*g);
    if z<>fail then
      if z=Size(veo) then b:=Size(veor)+1-y[z];
      else b:=y[z+1]-y[z];
      fi;
  l:=[1..a];m:=[1..b];
   
      for j in m do
        for k in l do
          
          if PositionSet(set,veor[y[x]+k-1]*veor[y[z]+j-1])=fail then nre:=nre+1;
          
          AddSet(set,veor[y[x]+k-1]*veor[y[z]+j-1]);
          set2[nre]:=veor[y[x]+k-1]*veor[y[z]+j-1]; 
          set3[nre]:=Concatenation(veoc[y[x]+k-1],veoc[y[z]+j-1]);
          fi;       
       od;
     od;
    fi;
od;
return [set2,set3];end;;
\end{verbatim}This code will be saved to file
makeset.txt
 \pagebreak

\subsection*{The function EqClassRel2(el)}
This function returns a set that contains the image under the
bijection $f:R_c\mapsto [1,88179840]$ of every element equivalent to
$el$ under relation $\approx$. The algorithm is evident and need not
to be explained.

\begin{verbatim}
EqClassRel2:=function(g) local h,set;
set:=[];	
for h in Mc do
   AddSet(set,PermToNr(h^-1*g*h));
   AddSet(set,PermToNr(h^-1*g^-1*h));
od;
return set; end;;
		
\end{verbatim} This code will be saved to file
eqclassrel2.txt
 \subsection*{The function
WriteRcv(set)}
 This function has as input a set of indices or numbers
and is called by command:
 \begin{verbatim}
gap>WriteRcv(set);
 \end{verbatim}
The command above takes every $n\in set$ and sets entry $n$ in the
"bit vector" to "false" or 0.
\begin{verbatim}
WriteRcv:=function(set)
local x;
for x in set do
 rcv[x]:=false;
od;
end;;
\end{verbatim} This code will be saved
to file writercv.txt

\subsection*{The function Redf(set)}
As mentioned earlier in the text the algorithm finds generator
expressions for elements in $E_f$ and each time it does this we need
to give it a $g\in R_e$. It was also mentioned that if given the
algorithm a $g\in R_e$ that is equivalent under relation $\approx$
with another $h\in R_e$ that has been before, then this $h\in R_e$
will not yield any new generator expressions. $Redf(set)$ returns a
set containing only one representative from each equivalence class
under relation $\approx$. The algorithm is simple and can be seen
directly by looking at the code.

\begin{verbatim}
Redf:=function(set) local x,y,set2,setemp,vec;
vec:=BlistList([1..Size(set)],[1..Size(set)]);set2:=[];
for x in [1..Size(set)] do
   setemp:=[];
   if vec[x] then 
      for y in Me do
        vec[PositionSet(set,y^-1*set[x]*y)]:=false;
        vec[PositionSet(set,y^-1*set[x]^-1*y)]:=false;
        AddSet(setemp,y^-1*set[x]*y);
        AddSet(setemp,y^-1*set[x]^-1*y);
      od;
      if IsEmpty(Intersection(setemp,set2)) then AddSet(set2,set[x]); fi;
   fi;
od;
return set2; end;;
\end{verbatim}This code will be saved to file redf.txt
 \pagebreak

\subsection*{The function CharToPerm(x)}
 Part of the $g\in R_e$ we are giving as input to the algorithm will
 come from the set $Redf(veven)$. With this set we find generator
 expression for 43984820 elements in $E_f$. To find expressions for
 the rest of 105100 we need some more $g\in R_e$ that has been
 extracted from a large list and are given at the end of this paper.
 This list took about 2 days to find on 20 computers. It consisted
 of representatives under relation $\approx$ of elements of the set $\{g\in
 R_e,L(g)=8\}$ and very few from $\{g\in R_e,L(g)=10\}$. From the
 last set only 17. Those were found using a solver found at
 \cite{re}. In order to compress the space the elements have been
 well coded and in order to make them permutations in $R_e$ we need
 the function $CharToPerm(x)$. If we want to convert an element of this list to a permutation we type in GAP:
 \begin{verbatim}
gap>CharToPerm("acekrihownvs");
(4,6,5,21,23,15,7,22,8,13,10,16,12,24)(14,20)
 \end{verbatim}
Since this is just something technical it needs no further
explanation.

\begin{verbatim}
CharToPerm:=function(g)
local s,v,w,m,k;
v:=[1..24];
v[1]:=11;v[2]:=17;v[3]:=19;v[4]:=22;v[5]:=13;v[6]:=8;v[7]:=24;v[8]:=6;v[9]:=18;
v[10]:=21;v[11]:=1;v[12]:=15;v[13]:=5;v[14]:=20;v[15]:=12;v[16]:=23;v[17]:=2;
v[18]:=9;v[19]:=3;v[20]:=14;v[21]:=10;v[22]:=4;v[23]:=16;v[24]:=7;
w:=[1..256];k:=InputTextString("abcdefghijklmnopqrstuvwx");
w[ReadByte(k)]:=1;w[ReadByte(k)]:=11;w[ReadByte(k)]:=2;w[ReadByte(k)]:=17;
w[ReadByte(k)]:=3;w[ReadByte(k)]:=19;w[ReadByte(k)]:=4;w[ReadByte(k)]:=22;
w[ReadByte(k)]:=5;w[ReadByte(k)]:=13;w[ReadByte(k)]:=6;w[ReadByte(k)]:=8;
w[ReadByte(k)]:=7;w[ReadByte(k)]:=24;w[ReadByte(k)]:=9;w[ReadByte(k)]:=18;
w[ReadByte(k)]:=10;w[ReadByte(k)]:=21;w[ReadByte(k)]:=12;w[ReadByte(k)]:=15;
w[ReadByte(k)]:=14;w[ReadByte(k)]:=20;w[ReadByte(k)]:=16;w[ReadByte(k)]:=23;
m:=InputTextString(g);s:=[1..24];
s[1]:=w[ReadByte(m)];s[v[1]]:=v[s[1]];s[2]:=w[ReadByte(m)];s[v[2]]:=v[s[2]];
s[3]:=w[ReadByte(m)];s[v[3]]:=v[s[3]];s[4]:=w[ReadByte(m)];s[v[4]]:=v[s[4]];
s[5]:=w[ReadByte(m)];s[v[5]]:=v[s[5]];s[6]:=w[ReadByte(m)];s[v[6]]:=v[s[6]];
s[7]:=w[ReadByte(m)];s[v[7]]:=v[s[7]];s[9]:=w[ReadByte(m)];s[v[9]]:=v[s[9]];
s[10]:=w[ReadByte(m)];s[v[10]]:=v[s[10]];s[12]:=w[ReadByte(m)];s[v[12]]:=v[s[12]];
s[14]:=w[ReadByte(m)];s[v[14]]:=v[s[14]];s[16]:=w[ReadByte(m)];s[v[16]]:=v[s[16]];
return PermList(s);end;;
\end{verbatim}This code will be saved to file
chartoperm.txt
\subsection*{The function MultEl(el)}
We write elements in generator form as a vector. For example [1,3,2]
means the element $UFL=s[1]*s[2]*s[3]$ and $s$ was defined above.
When we want a conversion from the form [1,..,.8] to a permutation
with cycles we use the function MultEl. Its code can be seen below:
\begin{verbatim}
MultEl:=function(el)
local x,an;
an:=();
for x in [1..Size(el)] do
an:=an*s[el[x]];
od;
return an;
end;;
\end{verbatim}

\subsection*{The function Invg(el)}
This function takes an element in generator form as defined above
and maps it to its inverse also given in generator form.
\begin{verbatim}
Invg:=function(g)
local set,x,y;
set:=[];y:=[];
y[1]:=7;y[2]:=8;y[3]:=9;y[4]:=10;y[5]:=11;y[6]:=12;
y[7]:=1;y[8]:=2;y[9]:=3;y[10]:=4;y[11]:=5;y[12]:=6;
for x in [1..Size(g)] do
set[Size(g)+1-x]:=y[g[x]];
od;
return(set);end;;



\end{verbatim}
\pagebreak

\subsection*{The main program}
The algorithm for the program was explained in the beginning and
here we give the implementation as GAP code:
\begin{verbatim}

sh:=0;cnt2:=1;positions:=[];
for x in sett do
    set1:=MakeSet(10,x);
    set2:=MakeSet(12,x^-1);
    cnt:=0;
    for y in [1..Size(set1[1])] do
        for z in [1..Size(set2[1])]  do
       
        g:=set1[1][y]*set2[1][z];
        if rcv[PermToNr(g)] then 
             positions[cnt2]:=Concatenation(set1[2][y],set2[2][z]);cnt2:=cnt2+1;
             set3:=EqClassRel2(g);
             WriteRcv(set3);
             cnt:=cnt+Size(set3);
        fi;
        g:=set1[1][y]^-1*set2[1][z]^-1;
        if rcv[PermToNr(g)] then
             positions[cnt2]:=
             Concatenation(Invg(set1[2][y]),Invg(set2[2][z]));cnt2:=cnt2+1;
             set3:=EqClassRel2(g);
             WriteRcv(set3);
             cnt:=cnt+Size(set3);
        fi;
        od;
    od;
    sh:=sh+cnt;   
    Print("Positions checked:",sh," ","Positions left:",88179840/2-sh,"\n");
od;
\end{verbatim}This code will be saved to file
mainprg2.txt
 \pagebreak
\section*{Running the program and output}
To run the program we make use of all above functions and sets. To
read those in memory and start the program the following code is
used:
\begin{verbatim}
Read("makesetv.txt"); #This file contains function MakeSetv
a:=MakeSetv(1);
vevenr:=a[1];
voddr:=a[2];
a:=MakeSetv(2);
veven:=a[1];
vodd:=a[2];
Read("makesetvc.txt");
a:=MakeSetvc();
vevenc:=a[1];
voddc:=a[2];
Read("makesety.txt"); #This file contains function MakeSety
y1:=MakeSety(1);
y2:=MakeSety(2);
Read("makeset.txt");#This file contains function MakeSet
Read("redf.txt");#This file contains function Redf
rcv:=BlistList([1..88179840],[1..88179840]);
Read("permtonr.txt");# This file contains function PermToNr
Read("eqclassrel2.txt");#This file contains function EqClassRel2
Read("writercv.txt");#This file contains function WriteRcv
sett:=Redf(veven);
Read("datav.txt");#This file contains the data included in the paper
Read("chartoperm.txt");#This file contains function CharToPerm
for x in datav do
AddSet(sett,CharToPerm(x));
od;
Read("inverseingens.txt");#This file contains function Invg
Read("mainprg2.txt");#This file contains the Main program
SaveWorkspace("test2");
quit;
\end{verbatim} This code is put in a file called
start.txt. And all the files(including datav.txt at the end of the
paper) in the $\backslash$gap4r4$\backslash$bin catalog. The program
is run under unix with command:
\begin{verbatim}
gap.sh -o 1000M start.txt
\end{verbatim}
The program will output:
 \begin{verbatim}
 Positions checked:3079007 Positions left:41010913
 Positions checked:18604155 Positions left:25485765
 Positions checked:19391883 Positions left:24698037
 Positions checked:22124037 Positions left:21965883
 ......
\end{verbatim}
And a couple of hours later...
\begin{verbatim}
 Positions checked:44089776 Positions left:144
 Positions checked:44089824 Positions left:96
 Positions checked:44089872 Positions left:48
 Positions checked:44089920 Positions left:0
 gap>
\end{verbatim}

\section*{What does this program do?}
As told in the beginning we need to compute each element in the
group $G_{cube}\cap R_c=E_f$ in at most 22 moves. This program does
this. When the program is finished the set "positions" contains
generator expressions of each element in $E_f$. An element of
"positions" has the form $[ 1, 2, 8,7]$. Here $1$ means $U$, $2$
means $F$ more general $n$ means $s[n]$. Actually not all elements
of $E_f$ are contained in positions but only representatives from
each equivalence class under relation $\approx$. To make a final
check that confirms that this set contains a representative from
each equivalence class under relation $\approx$ the following
program can be run:
\begin{verbatim}
Read("multel.txt");
corners:=Group(uc,lc,fc,bc,rc,dc);
rcv:=BlistList([1..Size(corners)],[1..Size(corners)]);
#rcv is a binary vector with all its entries set to true and size 88179840
for x in positions do
set:=EqClassRel2(MultEl(x));
WriteRcv(set);
od;
cnt:=0;
#the program below counts the number of entries in rcv that are true
for x in rcv do
if x then cnt:=cnt+1;fi;
od;
Print(cnt,"\n");
\end{verbatim}
 This program takes each element in "positions" converts it to a permutation with
disjoint cycles (this is the form which GAP can handle best). Then
build all elements in the equivalence class under relation $\approx$
and checks them of in the bit vector. After this procedure is done
the program counts the number of entries in the bit vector that has
been set to "false" and prints them out. If it returns 44089920 then
every element in $E_f$ has been checked of and we are finished. Note
that this final check actually proves the whole content of this
paper is correct. Actually this verification has already been done
by Bruce Norskog soon after i published the results of this work on
the internet.
\section*{Final comments}
It might seem that we were lucky that all elements in $E_f$ had
length 22 or less but the truth is that we had a backup plan if this
would be the case. We find generator expressions for each element
$g$ in the Rubik cube group in two steps. First $g$ is multiplied by
$b\in G_{cube}$ from the left and then $bg \in E_f$ and afterwards
$bg$ is multiplied by $a$ which takes it to the identity.
\[abg=id \Rightarrow g=b^{-1}a^{-1} \Rightarrow g^{-1}=ab\]
Assume that $a\in E_f$ has length 24. Assume that the diameter of
the Rubik cube group is 24. Compute the set \\
$X_1=\{aU^{-1},aL^{-1},aF^{-1},aB^{-1},aR^{-1},aD^{-1},aU,aL,aF,aB,aR,aD\}$.
Then each element of this set has length 23 or less. Then $a=ex_1$
where $e_1\in X_1,L(e_1)\leq 23$ and $x_1$ is a generator. Then
$b=x_2c$ where $x_2$ is some generator and $c_1$ is an element of
length 17 or less. When writing $g^{-1}=ab=e_1x_1x_2c_1$ we can
choose $x_1=x_2^{-1}$ and cancellation occurs giving that $g^{-1}$
has length 40 or less. It is well known that there exists a position
of length 26 so even this procedure might fail. So assume that $aU$
has
length 25 in the list above. Then compute the set \\
$X_2=\{aUU^{-1},aUL^{-1},aUF^{-1},aUB^{-1},aUR^{-1},aUD^{-1},aUU,aUL,\\
aUF,aUB,aUR,aUD\}$ and assume that all elements of this set have
length 24. Then $g^{-1}=ab=e_1x_1x_2c_1$. If $x_2\neq U$ we solve
this problem as before. If $x_2=U$ then
$g^{-1}=ab=\underbrace{e_2x_1U^{-1}}\underbrace{Ux_3c_2}$. And $c_2$
has length 16 or less. Again we can choose $x_1=x_3^{-1}$ and the
length of $g^{-1}$ is less than or equal to 40. Finally some
elements in the last set might have length 26 in which case we must
carry one last step of this. Obviously if there were too many
elements in $E_f$ having length 24 this computation could take
years. I must remark though that the elements of length 24 are very
rare and very few are known. In 1995-1998 the only known element of
length 24 was the "superflip" and today about 150 others are known
to the best of my knowledge.
\section*{Conclusions}
In this paper we have constructed a program that returns generator
expressions of length 22 or less for each element in the subgroup of
the cube group that fixes the edges. As discussed in the abstract
this means we have given a upper bound of 40 for the cube group.
\section*{Acknowledgements}
I want to thank my advisor Gert Almkvist for his help and
suggestions and my second advisor Victor Ufnarovski  for taking the
time to look through this work in detail. I further want to my
father, my wife and my mother for their help and support. Thanks to
Tomas Rokicki for publishing his result for the edges analysis and
Bruce Norskog who verified the output of this program. I also thank
Michael Reid for sharing his optimal solver. Finally thanks to the
people who contributed to GAP development, this program makes it
very easy to test mathematical ideas.
\section*{Data vector}
This should be placed in a file called datav.txt: \\
\\
\begin{tiny}
datav:=[ "aceglpmwqius", "acegpkwsujqm", "acejhtmoqluw",
"acejpkwmuqsg",
  "acekrihownvs", "acekvihqmsow", "acelismhqvpw", "acelismwohuq",
  "acelismwqhou", "aceliwqgonvt", "aceliwqugtmo", "acelqwmspgui",
  "acepikmsqhuw", "acerikmtohwv", "acesikwmqhpu", "acexigqkosmv",
  "acfgkojqmvxs", "acfkihmoqvtw", "acfunilqshow", "acgtikquwnoe",
  "achfikmoqwvt", "achliewuonqt", "achmikqxosfv", "achpiemsqkuw",
  "achpikfqmsvw", "achrvkmojsfx", "achseuqowkni", "achsikqownfv",
  "achtikuqonwf", "achuislmfxoq", "achvikoqfswn", "achvismoqekw",
  "achviswlonqe", "achxikmvprfs", "achxikufmsor", "acjfrnlxosvh",
  "acjgnrloewsu", "acjgnrlwesou", "acjgrnewlsou", "acjtlnhqevow",
  "acktieqowhvn", "ackxiegqmsov", "ackxiequosgn", "aclgniewprus",
  "aclgniuoexrs", "aclgwieqmsup", "aclheipqmvws", "acloniesqgwu",
  "aclonieswguq", "aclonisewgqu", "aclsvhmqjeow", "aclxrimohsfv",
  "acmgeujwoktq", "acmgkfjwprus", "acmgkfxqjsup", "acmgpkjwequs",
  "acmgsljwqvof", "acmoglqfwvjs", "acmokfjsqgwu", "acmokfjwqvtg",
  "acmteojwqkug", "acmxvgjoqetk", "acnfikhsqpwu", "acnfikhwqspu",
  "acnfikpxqsvh", "acnfikqhwspu", "acnfikqxovth", "acnkihvoqetw",
  "acolihmsweuq", "acoqihmewsul", "acoxikmvqsge", "acplietqmgux",
  "acplismwqveh", "acptieuqmkwg", "acpvikmqesxh", "acqgfnjwlsou",
  "acqokfmsjgwu", "acqsikuwmgoe", "acruikomxsge", "acsuikomgeqw",
  "acsvrihownke", "acsxikfoqngv", "actfikmqvhow", "actfikmuqgxp",
  "actfkmjwqhou", "actlvieqmhow", "actpieuqmkgw", "acugfnjwoslq",
  "acugkfojwsmq", "acuokfmswgjq", "acuokfsmjgwq", "acvfikmswphq",
  "acvfqjmkpsxh", "acvgpetkmjwq", "acvlhjmoqftw", "acvmietoqkxh",
  "acvmikqfosxh", "acvnkfxqjtgo", "acvriktomewh", "acvsilmoqexh",
  "acvwikmrosge", "acwtieqgokmv", "acwuiemoqslh", "acwuiksqmego",
  "acxfikhmupqs", "acxfikoqmsvh", "acxlisumfhoq", "adegmkiwpqus",
  "adegqkmupitx", "adeohkmwqius", "adesukhomiwq", "adeuokqsximg",
  "adeuqkmgtiow", "adewqomglius", "adexikmsqoug", "adfvslnohiwq",
  "adgfqkmvxitp", "adhuisrwpnle", "adkgmsfwpivq", "adkgqmewpius",
  "adkgwsfqmipv", "adknqgpuwifs", "adksowhfuiqn", "adktiomrhvfw",
  "admfkuhxrisp", "adpgikxeunqs", "adpgslmwrifv", "adpgslwfuirn",
  "adptkgnwqifv", "adrxsnlohifv", "adsvowqlhime", "advlmerwpgti",
  "advsiemkqoxh", "adxkmeuqpgti", "aegoidmlqsuw", "aehtikmcquwo",
  "aehtikmvqwdo", "aekgismvqwdo", "aelgniwcuqso", "aesoiwmvqclg",
  "aesuqkmgpwdi", "aevdikmhqwto", "aexoivtgmqlc", "afckiumoqgtw",
  "afkoismvqwgc", "afkrcsmxipug", "aflcqhpuwims", "aflcqhupmiws",
  "afqckhjuwpms", "agclnpewqjus", "agcuikmoqswe", "agdsiemoqluw",
  "agdsikfqmxov", "agecltqowimu", "agecqkmpxius", "agecskmoqiuw",
  "ageisdmoqluw", "agejckmwqpus", "agejktmoqcuw", "agekciurmpws",
  "agelcjqosumw", "agelcjurmspw", "agelcpmrixus", "agelipquwcms",
  "ageliqmwocut", "agelismoqvdx", "agenckuoitwq", "agenckuqitow",
  "agenikqxocvt", "agenikuwqcot", "agenksrdpiuw", "ageoikmdqtuw",
  "ageoikmrxsuc", "ageoikmsqwuc", "ageoikmsqxud", "ageoikmsxcur",
  "ageoitmlqcuw", "ageonkqcjsuw", "agepltmricuw", "ageqckuniswp",
  "ageqxkmcosuj", "agesckuomxri", "agesikmcqwuo", "agesikmrxcuo",
  "agesldmoqiuw", "agesnkuwocjq", "agesqkmwcoui", "agesqkpuwidm",
  "agesqkuomidx", "agesrkjomcuw", "agesrkjwmcou", "agesvkjwocmq",
  "agesvkojmcqw", "agesxkmoucqj", "ageuckiosxmr", "ageuckqoixmt",
  "ageuikmwqsoc", "agevismoqlxd", "agewikmsoquc", "agewikmsqvpc",
  "agewikmvqsdo", "agewikqucsmo", "agewikucmsoq", "agewikuqmtdo",
  "agewmkuordsi", "agewpkmjcsuq", "agexikmcpsuq", "agexikmcqpus",
  "agfkvqxdmitp", "agflniquoctw", "agfsnivkqcpw", "agjsfvmwoclq",
  "agknwsjoqcfv", "agkrcsmoivfx", "agksiqfomcwv", "agksiwmqecvp",
  "agksrwmfocjv", "agktipmwqcfv", "agkuiwtoqcne", "agkvpsfncqwi",
  "agkvqsmoxifd", "agkvsfjoqcnw", "agkvwsmpqidf", "agkwiomcqsfv",
  "agkwiomrtcfv", "agkwiomsqfdu", "agkwismcqvfo", "agkxcjmoqsfv",
  "aglcrtmoweju", "aglfixuomctq", "aglovieqmscw", "agloximcquse",
  "aglrniewpcus", "aglsjweoqcmu", "aglwvieomqsc", "agmctkjoqeuw",
  "agmldpjwqeus", "agmodejvqktw", "agmsikqoecuw", "agmsikqwecou",
  "agmsikuwoceq", "agmwkfjcqsuo", "agnfiwquoctl", "agnrsljowceu",
  "agodikmvqswe", "agplismwqcfv", "agpsikfwqcmu", "agqsikeomcuw",
  "agqsikewmcou", "agscikmoqeuw", "agsjrkeowcmu", "agswnielqcuo",
  "agtdikmwprfv", "agtovqedmikx", "aguerisowcml", "aguliemoqcxs",
  "agusikewocmq", "agusikfdpnwq", "aguwkfjqmsoc", "agvcelmoqjtw",
  "agvcflmoqitw", "agvcfpmiqktw", "agvfikmcqspx", "agvfkmjoqctw",
  "agvkiemcqsxo", "agvkiemwqods", "agvrcemoiltw", "agvrcemoitxk",
  "agvsipmxqcle", "agvwiemcqkto", "agwsipmlqcue", "agwvieqtocmk",
  "agxkcetnupqi", "ahecitmrkpuw", "aheodirsxnul", "ahesqkmowdui",
  "ahnuieqcowlt", "aidkmersowhv", "aiegckqwnpus", "aiegquptwcmk",
  "aiehmutqpldx", "aiekqumdpsxh", "aielqspdwnvh", "aielqspudhmw",
  "aienckwguqsp", "aienpkwcuqsg", "aieomkugqwsc", "aierckmsgwup",
  "aierpkmscwug", "aieskhmdpruw", "aiesmkuwqgoc", "aiespkmdgquw",
  "aievqkpdwhms", "aiexckmgqusp", "aikgcsmqvpfx", "aikhscqoxnfu",
  "ailgncesqxup", "aimgekupqsdw", "aimgslpfdvrw", "aiswceukmroh",
  "aiurewdkmptg", "aiurkfdgmwsp", "aivgcewktnqp", "aivrcehkmwsp",
  "aixfmkptdhuq", "ajchkfosxqnv", "ajegpkuwcmsq", "ajegxkqtcpmu",
  "ajeovkmwcqsg", "ajesqkpdhmuw", "ajewtkucmqgo", "ajkgxsmrcpfv",
  "ajkvmsfdhqwp", "ajndwepkvqtg", "ajogmewcuqtk", "ajsgpvqlcwme",
  "akecqhiwmpus", "akecwhirosmu", "akehctmoiquw", "akeoiguqmwsc",
  "akeoiqncwshu", "akerchmsixup", "akeritmwphuc", "akeschmriouw",
  "akesiwugmcoq", "akhpcsmwiqfv", "akhtqmewpjdv", "akmchfjrosuw",
  "akmuetjcqwph", "akqhcemsiwup", "akvcmetiqpwh", "akvncetgiqwp",
  "akvscemriphw", "akvsiehqmcow", "akvtweqgpimd", "akvwceitgpmr",
  "akxoivtemqhc", "akxqcemgipus", "alcgifmwqous", "alcoifugmqsw",
  "alcsrfmowgju", "alcsvfmojgwq", "alcsvfmqjgow", "alegctirwpmu",
  "alegmtudpiwq", "aleqptmjgcuw", "alesdiqowgmu", "alesdiuomgwq",
  "alesixmoqcug", "alesoidrnguw", "alhsiemoqcuw", "almehtrdpiuw",
  "almgqsjwpdfv", "almsfdjoqguw", "alpscfmriguw", "alsgivmwqoec",
  "alsgrieowcmu", "alsgvieomcwq", "alvntfjoqcxh", "alvsdemoqgxi",
  "amcgikouwqes", "amegiksqcuwo", "amegikuqcwso", "ameoikwruchs",
  "amepqkwguisd", "amesikgocquw", "amesikocwuqg", "amesikoqugwc",
  "ameuhkscjqwo", "amlgipewqcus", "amrfikuhscwp", "amswpvceqjlg",
  "amugkfwcjqso", "andtgvlqipfx", "anegikqcwups", "anigkfquwpds",
  "anksiwqhdufo", "ankwisfgqcuo", "anlgriewpcus", "anlpchuqiews",
  "anlvgqepxsdi", "anpuiglcqsfw", "anvkjcewqpth", "anvwcepqitgk",
  "aoclrfwjugms", "aocuikqswgme", "aoecikquwtmh", "aoecikuqmtwh",
  "aoecktmwqguj", "aoednskwqjuh", "aoedqkmvhitw", "aoedukwirgsm",
  "aoegdmjwqkus", "aoegikmrdvxs", "aoegikqsucmw", "aoegikqtvcmw",
  "aoegikquwntc", "aoegiktqucwn", "aoegikuxmctq", "aoegikwqmctv",
  "aoegisukmcwq", "aoegltmwqiuc", "aoegnkjqucws", "aoegskmqjcux",
  "aoehskmdqiuw", "aoejmkdwqgus", "aoejrksuwcgm", "aoelisdqmhwv",
  "aoelismdqhxu", "aoelituqmcwg", "aoelnrcuwgjs", "aoelqwmvhitc",
  "aoenkgjvqctw", "aoenskwjvcqh", "aoerckmwigut", "aoerckmwjgus",
  "aoerckwmuhis", "aoerclmwigus", "aoerjkdwngus", "aoerskmxicug",
  "aoesihmwqkud", "aoesihquwcml", "aoesikmcquwg", "aoesukwrhcim",
  "aoetikmgqcwv", "aoetikqhwcmu", "aoevhkjwqcns", "aofgimuwqctk",
  "aofgksquwcni", "aoflsdjwqgnv", "aofriktwmcuh", "aogcrnlwfius",
  "aogmvilqecws", "aogrldewnius", "aohckfjwqnus", "aohclsqfwimv",
  "aohlivmwqcfs", "aohlnevwqjds", "aohlvemqjcws", "aohlweiuqcms",
  "aohtifmwqcul", "aohtikfdmvrw", "aohwikmuqcfs", "aoignsquwckf",
  "aoigrkewncus", "aoilnevwqchs", "aojgfnlsxruc", "aojgkvcqmews",
  "aojglrwfuctn", "aojgnveqlcws", "aojkngewqcus", "aojtndkwqgfv",
  "aojwfnqcguls", "aojwgvsqmcke", "aokdqshfwimv", "aokgcsirvxmf",
  "aokgismvwcfq", "aokgiswufcmr", "aokgixmuqcfs", "aokgmsqwicfv",
  "aokgrsjwmcfv", "aokgvsjfqcmx", "aokgxsqfucjn", "aokjdemwqgus",
  "aoknigfwqcus", "aokrsgjfwcmv", "aoksivmwqcfg", "aokuismhqcfw",
  "aokvishfrcmw", "aokwiemgqctv", "aokwismgqufd", "aolcniewqtuh",
  "aolgnvjwqces", "aolgwieqncus", "aolhismwqcfv", "aolmeiswqcgv",
  "aolsnihurcew", "aoltiwhurcfn", "aoluiqswncge", "aomchvlqejws",
  "aomckfhwqius", "aomckfjwqtuh", "aomedkgwqjus", "aomedkjwqtuh",
  "aomgikeuwcqs", "aomgxtjlqcue", "aomidfhwqlus", "aomligrwecus",
  "aomtfljwqgdv", "aomtilcwqgfv", "aonfikdsqhwu", "aonfikqxhcvt",
  "aonfwkqtucih", "aonlisqfwgdu", "aonrceiuwhtk", "aonuiewckgqt",
  "aonuisfwqcgk", "aonxieqkvctg", "aoqcehjuwkms", "aoqgntwlucje",
  "aoqhniesxckv", "aoqhskjfwcmv", "aoquismlwche", "aosgrvmlwcje",
  "aosgvxmqjcle", "aotfikmdwhuq", "aotliemdqhuw", "aotungkwqcie",
  "aougfrljmcws", "aougriwlfcsn", "aoulcfnqjgws", "aounkfdqihws",
  "aoutigqlwcme", "aovcgemwqitk", "aovciemwqshk", "aovcikmhqsxf",
  "aovfikmhqcwt", "aovfikmwqtgd", "aovgemjdqktw", "aovliemtqcxg",
  "aovrceitwhmk", "aovsikmdqgxf", "aovtikmgqcfx", "aowckhmeujqs",
  "aowgnielucqs", "aowjkfmuqgsc", "aowliemsqcgv", "aowuieqtglmc",
  "aoxcikfqmtuh", "aoxgniesqcuk", "aoxrceuimhtk", "apegikqwsnuc",
  "apegmkuqxisc", "apejmkugxqsd", "apelsdmriguw", "apevckqxihms",
  "apknisfvqcwg", "apvnkfdgiwqt", "aqcwikmfsogv", "aqegcomixlus",
  "aqegikmdsowv", "aqegvkmcjwso", "aqelismdxouh", "aqencksuwpig",
  "aqeoikmgcsuw", "aqeoikutmwdg", "aqesvkmwjgoc", "aqevikxsmodh",
  "aqewikucmosg", "aqewvkmgjcso", "aqjsfnuwlgoc", "aqkgosmiwufd",
  "aqkncgutixfo", "aqlwniuoesgc", "aqnfpkcjwsvh", "aqtvkfjgmwdp",
  "arcgienuwptk", "aregxkumdpis", "areldinuwphs", "arflgsmwipdv",
  "arhojckusxfn", "arlgcinuwpes", "artocvhimxle", "arvwceipgtmk",
  "asedikupmqwg", "asegikwpudnr", "asegipulmqwc", "asegixqulcmo",
  "asegvkmqjwco", "asekigmoqcuw", "asekiumgqwdp", "aselcwmgipqu",
  "aselihmoqcuw", "asenckigqxuo", "aseoikuqgmwc", "aseoikuxmcgr",
  "aseqikmogcuw", "aseqipulmcwg", "aseqpkujhmwd", "aserckmviwog",
  "aseuikqcwomg", "asevikmoqchw", "asewcqmlipug", "asewikmogquc",
  "ashxikfoqncu", "askoigfwqunc", "aslcngeiqwuo", "aslorimweugc",
  "aslqvimwgceo", "aslrniegxouc", "aslwipmgqeuc", "aslwnieohcur",
  "asnuiepwqcgk", "asqofnlugwjc", "asuwikfqmcpg", "asvkiemgqwpc",
  "asvwiemkqodg", "aswvioqlgcme", "asxhikmcqofv", "ategikmvqowc",
  "atelcuqhiomx", "ateliwmchruo", "ateqhkmociuw", "aterxkigcvmo",
  "atngceikuxqp", "atnuweqkpdjh", "atpgikudmqfw", "atvgpeikcxmq",
  "atvwiemkqodh", "auegckqswnip", "auegikmxqotc", "auegiknswrdp",
  "auejpkhwdqms", "aueoikqgtxmc", "aueoqkmsxdgi", "aueoskigqwmc",
  "auesikmwqocg", "auesrkiwhomc", "aufmriskwodh", "auhtwkinqofc",
  "aulsniqogcew", "avcfikmsqwph", "avegckmsiwor", "avegikmcqspw",
  "avegikmsqodw", "avegikomqsdw", "avegikqdswmo", "aveligmoqsdx",
  "aveoikmcqsxg", "avesikmwqodg", "avesqkmdpwjg", "avkocgfmixqt",
  "avkpismhqwfc", "avkwismohcfr", "avkwmsfdhioq", "avwsiemgqcpk",
  "awdjlhmsqofv", "awedismkqouh", "awedqgpukims", "awedqkphvims",
  "awegikmdqosv", "awegikmuqodt", "awegikomsuqc", "awegikqscumo",
  "awegmkurpcsi", "awegqkmipvtd", "awegskqupdmi", "awejckosnqug",
  "awekisuqmodh", "awelivmhqods", "awelpdmsirug", "awenckvsiroh",
  "aweoihmsqluc", "awerckusmioh", "awesikuomgqc", "aweuikmoqsgc",
  "aweuikmsqcgo", "awgjpvmscqle", "awjgfvmsqolc", "awjrcesmuphk",
  "awkoqsmduifg", "awkoqsmvhifd", "awmfikhvqods", "awmoekjdqgus",
  "awmokfjgqusc", "awmuktrepdgi", "awpfikmuqcth", "awrfikdusomh",
  "awsehvmlqcjo", "awsmikgeqouc", "awsoeqjmuchk", "awvoremtkcjg",
  "awvsikfqmodh", "axegikusmqdo", "axeoqkmsgdui", "axhtikqfucpn",
  "axksiqocfumg", "axphckifsvqn", "axqgemitupdk", "bdmuqtlxjoeg",
  "bfmsokjhqcuw", "bkeogjmdsxuq", "bwehvqkpjdnt", "caljnpegqsuw",
  "camlxfroigus", "cbsgivlqmpwe", "cbvfqjmwplgs", "cealnigousqw",
  "ceglibmvqwps", "cejlnpawqgus", "cesagjuwmqlo", "ceslibmguwqo",
  "cfsulranwoig", "cfsuqjmbphxk", "cgalfiuoswmq", "cgaliqewosun",
  "cgalpiuemsqw", "cgalqiemuswp", "cgaltiqouemw", "cgaltiuomeqw",
  "cgameiwoqslv", "cgateiuqknow", "cgauniswoelq", "cgaurimwseol",
  "cgbxkimoqsfv", "cgeiqpmbxlus", "cgelibmvprtw", "cgelibqtowmv",
  "cgelibvqmsox", "cgelibwusqom", "cgelmjpwbsuq", "cgelniuwosaq",
  "cgelpbujmsqw", "cgelqbouismw", "cgelwbiornus", "cgesqamokiuw",
  "cgetibmlovwq", "cgfbniqosuxk", "cgfxblniqsov", "cgfxikmoqsbv",
  "cgjlbfqwmsou", "cgjlfsaoqxnv", "cgkasjqowvfn", "cgktibvmfqow",
  "cgmakpjwqeus", "cgmesbjoqlxu", "cgmilfrotauw", "cgmpbljrtwfv",
  "cgmxufaoqkis", "cgnlibfowrvs", "cgnlibfqusow", "cgofibqmxsvk",
  "cgpilnfqavxs", "cgplibufmsqw", "cgplibuqfsmw", "cgpsaeuqmixk",
  "cgpxikfqmsub", "cgrfibwuokmt", "cgsfibmuqkow", "cgsubfjwomlq",
  "cgsuibqfomkw", "cgtbniuoeqxk", "cgulibfnqsow", "cgvbiemoqlxs",
  "cgvlbfrojnxs", "cgvnafrokixt", "cgvsiembqokw", "cgvtibwoqkfm",
  "cgvxniqboetl", "cgwlbfqjosmu", "cgwlibeqosum", "cgwlibqeusom",
  "cgwlifmoqvta", "cgwtblqfjvom", "cgxfnbutokjq", "cgxtibmoqlfv",
  "chelibpuwnrs", "chesikunrwob", "chvnqfptwbjk", "chxlqemtpiub",
  "cimalfrwpgus", "cixkemugpqtb", "ckafnisgquwo", "ckafsiomvwgq",
  "ckahniuseqwo", "ckgbniseqwpv", "ckwhibmoqvte", "ckxpieuqmath",
  "claovimgsweq", "clebiomgqsuw", "clergbmoisux", "clewgbmrisup",
  "clewibmogsuq", "clmgbfusjwoq", "clmsbfjwqgou", "clqajfhuosmw",
  "clqsbfmwjgou", "cluwbfjqmsgo", "cmalriseguwo", "cmulbfjqgwso",
  "cnebiqkpvgws", "coafnlrwigus", "coafsjmwqguk", "coaheikmuxqs",
  "coalnwiuqges", "coalriesugmw", "cobgilqfwsmv", "coegikuamqws",
  "coegmkuaqiws", "coehibmuqxtk", "coehibukrnws", "coehqbkuwims",
  "coeiahmwlqus", "coelrbmsvxjg", "coelsbujmgqx", "coergaiuwlms",
  "coernbtjvgwk", "coesibmvqlhw", "coetbgjwqlnv", "coetkbxqjgnv",
  "coglbtjqmuwe", "cohtbljqmuex", "cohuibfqmwtk", "coifrbvwngtk",
  "cojhenaqkvws", "cokfibnqsgwv", "cokithewnrbv", "cokjnaquwges",
  "cokrvbmxjgfs", "coktibfhuxnq", "coktibfqmxug", "coljnberhsuw",
  "comhbljfwtuq", "comlfubqjgws", "conbrixutgle", "confaktwqiuh",
  "confibtkquwg", "confibuqwksh", "conwibqfugtk", "coqlbfwjmgsu",
  "coqlfnsjwgau", "coqlibesmgwu", "coqlnbvwjgfs", "coqlviweagms",
  "coqtbfjukgmw", "coqwbtjufgln", "corlibwmfvtg", "corlibwusneg",
  "cortliawfgnv", "cosfibqmwguk", "cosnbfjuqgwk", "cosqvnweaglj",
  "cosubkmwqfgi", "coulbfsjqgmw", "coulibwmqges", "covfibhwnktq",
  "covfmbhwrkti", "covfnbkwqgjt", "covgnaewqski", "covlrbjwegns",
  "covwisfqmglb", "coxfibumtgqk", "coxknbteugqi", "coxtiulqmgbe",
  "cqaoxiemugsl", "cqaxnikovgfs", "cqehmbwpuisl", "cqelibuspgwn",
  "cqewihusmola", "cqjbenxskoug", "cqmlsbvxjofg", "cqsoibmlwgue",
  "csagniuoeqwl", "csahnqelpwui", "csalniouwqge", "csaqnikohwfv",
  "csawjkegqnuo", "csfbnipgqwuk", "csnhibeowkru", "csqlubigwemo",
  "csvfrbmxgkjo", "csvmikqgwfao", "ctbrgoixlnue", "ctelawmiqoug",
  "ctvgweikqbmp", "cuakripswgfn", "cueagsmpirxl", "cuelibmgxqto",
  "cukgibqnosfx", "cunfaoiswrhk", "cvbrgsmliwfp", "cvksibmoqwfh",
  "cvktfiqgmwbo", "cvmpbtlqjwge", "cwehbujsqkno", "cwelmbusjgor",
  "cwelqbumiosg", "cwevnbikqotg", "cwjlnveqaosg", "cwkhbsjqmoev",
  "cwluiomsqhae", "cwpakimsqhfv", "cwqpbljfsvmg", "cwvhifqbkoms",
  "cwvnbfjrpgtk", "cwvnofrsjlhb", "dbegikuqmows", "ebmgksrwpiud",
  "ecmlsajoqguw", "ecpnbkwjusrh", "ecpnbkxjusqh", "ecsxbknqjvgp",
  "edmkqbhwpius", "edphknasxiqv", "egloibmcqwsu", "egsaidwmouql",
  "egsvcojrbnlw", "egvkbmitrpdw", "ehpkqbjsxdnv", "ejulmbrpgdws",
  "ejusbkgocmwq", "ejusbkndhqwp", "eknscqjhbovx", "ekvgbmcqwtoj",
  "elaoidugmqsw", "enkobgqvswjc", "eogcxkumjaqs", "eokgcsjrbxnv",
  "eomgjkdwqaus", "epksbwduihrn", "epvgbmtkrdwi", "erhucsblwpjm",
  "erxgbvtidplm", "esuwbkmogcjq", "euigbkxsqomc", "euinbkqgswdp",
  "euinbkqtogxc", "evajpdwnglqs", "ewctmibvkqoh", "ewmockjhbrus",
  "fcagsqjolnuw", "fhspbmdviwlr", "fhvcqjwmtapk", "fonaidqkwshu",
  "gcajntuqelow", "gcalitmoqeuw", "gclantuoejwq", "gdeaqjmswlup",
  "gemxdaqvoljt", "gfurbdjomksx", "gipkecmqbtwv", "gkanqdewpjus",
  "gkcanqpuwies", "gkeaspquximc", "glcaxiouesmq", "glebmurdpsxi",
  "glmrdftojbux", "glnfidqbwspu", "glpxvcbimsfq", "gnparceuwlis",
  "gobwidkqmtfu", "gockrbjwmtfv", "gojsdvaqmwle", "goqjnbxudles",
  "gosurdmwelja", "govlidmwqbtf", "gpeiqdusmbwl", "gquednaljwso",
  "gsajvdolewqm", "gseiwcqapuml", "gtdkqjmwpafv", "gtnuidqeolxb",
  "guitdbqnwofl", "guktidfqmbow", "gwbtmjpvqdfl", "gwecipquslma",
  "gwvfidmlqsao", "gxplrcnfitub", "hcofibqmusxk", "hrajsdfnkxov",
  "ibfogqvxmcsk", "ibhoqkmvxcft", "ibpgqkfdmvxs", "ibptkgnwqcfv",
  "icbhqwmsplfv", "icelqwmgosub", "icvomehwrksb", "idaukmsqheow",
  "idegqkusmowb", "idetqlwmbgpv", "idkhbfrnxout", "idsulrafpgwn",
  "idvoqextbgml", "idvtqblwpgfn", "iektqbmvpwhd", "iekwbvgscorn",
  "ifvkqbmsxodh", "igqnblpfdsuw", "ihbksdwfvoqn", "ihclnqpuwfbs",
  "ihpbnqukedxs", "iksuqemhxbdo", "ilbxqdhuofms", "ilvbqemdxohs",
  "imrfwbuhsdpk", "iodglbmwqtfv", "iomfbutxrgdk", "ipelgbusqnwc",
  "irsbdhlvmeow", "isalnqdguexo", "isepqbxukdmg", "iuehqbwmkdps",
  "ivatnqhfdwko", "ivbpnqewkgdt", "iwtfqbmkpuhd", "jnakwuqdsgpe",
  "joqgebvwmctk", "jragvkpwmces", "jskwnqugecbp", "jsowndegqlua",
  "juaonetwqchk", "kdjebnrohsuw", "kdqgwsanpfjv", "kesugdinoqaw",
  "kgiebdmoqvtw", "khmaqdjofsuw", "kibtgfnqdvow", "kixfqbupmcth",
  "krmxcsjobgfv", "ksmwcuafhqpj", "lceoshiaqumw", "lecunipaqshw",
  "licmgovbfwtr", "lnadqtejhwuo", "lphvcsjrbenw", "lwepchsmuria",
  "lwjbhfrspdnv", "majedvoqglws", "mgarcoelisuw", "mhqlvawtepjc",
  "mjeskbuqcgow", "mjpfrbukcgxs", "mleapdugsqwj", "mocwblvqjhte",
  "mougkiawqces", "mqoclsfjavgw", "mruagfjpdlws", "mrugfkwjdpbs",
  "ngibqdxlsufo", "obxngvtciqkf", "ojdkxrmsgafv", "olfbndquwgti",
  "omriwdtbkvgf", "pdmlwfargsju", "pgbicnqxfsvk", "prjncuawfkhs",
  "psmwcfjlbgur", "pvkacimgswfr", "qameowucjlsg", "qjecwsmgoaul",
  "qkenixvocsha", "qwvfcbmoikhs", "rnoticvelxag", "scegaomrixuk",
  "sgelcbqowimu", "sgwoniuqelac", "tgeawdmoilur", "trgiwpmbeclv",
  "twvoiemhqadk", "txkpbfdgiurn", "unkgasfqiwco", "utbwiemhqocl",
  "vjtkoqmbhcex", "voekbsjwqcng", "wjeovkmaqcsg", "wseclpmgaruj",
  "xoahcqesniul", "xsvnafjgkrdp" ];
\end{tiny}

The last data have been selected from a list that would take about 2
months to generate on an usual PC but thanks to my advisor Gert
Almkvist access to fast computers have been given and made it
possible to compute in just a couple of days.

\end{document}